\documentclass[oneside]{article}
\usepackage[utf8]{inputenc}
\usepackage{amsmath, graphicx}
\allowdisplaybreaks
\newtheorem{theorem}{Theorem}
\title{Differential Equations: A Historical Refresher\footnote{This is a preliminary version of a talk delivered in an ISTE sponsored  Faculty Development Programme on Applications of Advanced Mathematics in Engineering held in Vidya Academy of Science \& Technology during 27.11.2017 - 01.12.2017.}}
\author{{\bf V. N. Krishnachandran} \\ Vidya Academy of Science \& Technology, Thrissur-680501 \\ email: {\tt krishnachandranvn@gmail.com}}
\date{}
\begin{document}
\maketitle
\tableofcontents
\newpage
\begin{abstract}
This paper presents a  brief account of the important milestones in the historical development of the theory of differential equations. The paper begins with a discussion on the date of birth of differential equations and then touches upon Newton's approach to differential equations. Then the development of the various methods for solving the first order differential equations and the second order linear differential equations are discussed. The paper concludes with a brief mention of the series solutions of differential equations and the qualitative study of differential equations. 
\end{abstract}
\section{Introduction}
Neither while learning differential equations at college nor 
during my initial years of teaching differential equations in
 colleges was I excited about this particular branch of mathematics. The theory looked like a bag of tricks mysteriously
 producing the correct answers always, especially the ``$D$
 thing''! It was then I stumbled upon cheap Indian reprints of two classic books on differential equations which altered my perception of this area of mathematics. Much of contents of the book by George F. Simmons\cite{Simmons} could be read (without pencil and paper!) with much pleasure, enjoyment and excitement and it gave a lot of insight into the minds of the early pioneers in the area. The book by Coddington and Levinson \cite{Coddington}, though written in a terse style, clearly and emphatically articulated the fact that the theory of differential equations is not about a bag of tricks, but there is a deep conceptual framework behind it. 
These notes are intended  to convey some of these notions to mathematics teachers, those who apply mathematics in their areas of specialisation and to all those who are bored with the current pedagogical approach to teaching differential equations.

\section{The date of birth of differential equations}
Some historians of mathematics  consider the day on 
which 
Leibniz\footnote{Gottfried Wilhelm  Leibniz (1646 - 1716) was a German polymath and philosopher who occupies a prominent place in the history of mathematics and the history of philosophy, having developed differential and integral calculus independently of Isaac Newton. Leibniz's notation has been widely used ever since it was published.} solved and wrote down the solution of the following, then ``remarkable" but now trivial, differential equation as the day on which the theory differential  equations was born:
$$
\frac{dy}{dx}=x.
$$
This happened on 11 November 1675 and this date is considered as the date of birth of the theory of differential equations (see \cite{Ince} p.537 and \cite{Sasser}). Perhaps what is more important is the way Leibniz recorded his solution as
$$
y = \int x\, dx = \frac{1}{2}x^2.
$$
It was the invention and use of the integral sign ``$\int$'' that made it all the more memorable and useful.
\section{Some earlier ideas: Newton's approach}
\subsection{Newton's classification of differential equations}
The search for general methods for solving differential equations had actually begun a little earlier by Newton\footnote{Sir Isaac Newton (1642 - 1726) was an English mathematician, astronomer, theologian and physicist who is widely recognised as one of the most influential scientists of all time and a key figure in the scientific revolution. His book Philosophiæ Naturalis Principia Mathematica ("Mathematical Principles of Natural Philosophy"), first published in 1687, laid the foundations of classical mechanics. He shares credit with Gottfried Wilhelm Leibniz for developing the infinitesimal calculus.}. Newton classified differential equations into three classes as follows:
\begin{enumerate}
\item
$\dfrac{dy}{dx}=f(x)$ or $\dfrac{dy}{dx}= f(y)$ 
\item
$\dfrac{dy}{dx}=f(x,y)$
\item
$x\dfrac{\partial u}{\partial x} + y\dfrac{\partial u}{\partial y} =u$.
\end{enumerate}
Newton, of course, did not write down these equations in these forms. He presented the equations in his notations of fluxions. These ideas are contained in {\em Methodus fluxionum et serierum infinitarum}, written around 1671, but published posthumously
in 1736 only. An English translation of the book is available online (see \cite{Newton}). In this work, Newton solved several first order differential equations two of which are given below (see \cite{Phaser}):
\begin{enumerate}
\item
$\frac{dy}{dx} = 1-3x+y+x^2+xy $ with the initial conditions $y(0)=0$,  $y(0)=1$ and $y(0)=a$
\item
$\frac{dy}{dx} = 3y-2x+\frac{x}{y}-\frac{2y}{x^2}$ with an indication of the general method.
\end{enumerate}
\subsection{Newton's method of solution}
Newton's general method was to develop the right-hand member 
of the equation in powers of the variables and to assume as a
solution an infinite series whose coefficients
were to be determined in succession.
To illustrate the method we examine how Newton solved the first of the above two differential equations:
\begin{equation}\label{eq5}
\frac{dy}{dx} = 1-3x+y+x^2+xy, \quad y(0)=0.
\end{equation}
As an initial approximation we take
$$
y=0.
$$
Substituting this in Eq.\eqref{eq5} and retaining only the lowest degree terms we get
$$
y^\prime = 1.
$$
This gives
$$ 
y = x.
$$
Substituting this in Eq.\eqref{eq5}, and again retaining only the lowest degree terms we get
$$
y^\prime = 1-2x
$$
This gives
$$
y= x-x^2.
$$
repeating the procedure we have
$$
y^\prime = 1-2x+x^2
$$
and hence
$$
y=x-x^2+\frac{1}{3}x^3.
$$
Continuing this way, Newton obtained the solution as
$$
y=x - x^2 + (1/3)x^3 - (1/6)x^4 + (1/30)x^5 - (1/45)x^6 +\cdots
$$
It may be interesting to see the closed form solution of the equation:
$$
y= a\left( \text{erf}\,(bx+b)-\text{erf}\,( b) -\frac{4}{a}\right)e^{\frac{1}{2}x(x+2)} + 4 -x
$$
where 
$$
a= 3\sqrt{2\pi e}, \quad b= \frac{\sqrt{2}}{2}.
$$
Note that this involves the error function which is not an elementary function.

\section{First order differential equations}
The various general methods for solving differential equations were developed as responses to the challenges of solving practical problems. 
\subsection{Variables separable type}
\subsubsection{Formation of the concept}
Any differential equation which can be put in the form
$$
f(x)\, dx = g(y)\, dy
$$
where $f(x)$ is a function of $x$ only and $g(y)$ is a function of $y$ only is said to be of the variables type. The solution of the equation is given by
$$
\int  f(x) \, dx = \int g(y) \, dy.
$$

Now every schoolboy knows this and it is as clear as sunlight! But not so for mathematicians of  the seventeenth century.

A differential equation of the variables separable type arose for the first time in the solution of the problem of isochrone published by James Bernoulli\footnote{James Bernoulli (also known as Jacob Bernoulli (1655  – 1705) was one of the many prominent mathematicians in the Bernoulli family. He was an early proponent of Leibnizian calculus and had sided with Leibniz during the Leibniz–Newton calculus controversy. He is known for his numerous contributions to calculus, and along with his brother Johann, was one of the founders of the calculus of variations. He also discovered the fundamental mathematical constant $e$.} in 1690. The problem is to find the curve in a vertical plane  along which a body will fall with uniform vertical velocity. The curve is called the isochronous curve. James Bernoulli's method required the solution of a differential equation of the following form:
$$
\frac{dy}{dx} = \sqrt{\frac{a}{by-a}}.
$$

It was Leibniz who discovered the underlying principle of separation of variables and he communicated his discovery as a great theorem to Huygens\footnote{Christiaan Huygens (1629 – 1695) was a prominent Dutch mathematician and scientist. He is known particularly as an astronomer, physicist, probabilist and horologist. 
His work included early telescopic studies of the rings of Saturn and the discovery of its moon Titan, the invention of the pendulum clock and other investigations in timekeeping. He published major studies of mechanics and optics (having been one of the most influential proponents of the wave theory of light), and pioneered work on games of chance.} towards the end of the year 1691. To John Bernoulli\footnote{ John Bernoulli (also known as Johann Bernoulli) (1667 – 1748), a younger brother of James Bernoulli, was a Swiss mathematician and was one of the many prominent mathematicians in the Bernoulli family. He is known for his contributions to infinitesimal calculus.} is due the term and the
explicit process of {\em seperatio indeterminatarum} or separation of variables.

But it was not smooth sailing always.  It
was noticed that in one particular yet important case this process broke down; for
although the variables in the equation
$$
axdy - ydx=0
$$
are separable, yet the equation could not be integrated by this particular method.
The reason was that the differential $dx/x$ had not at that time been integrated!

The then newly discovered tools of calculus were put to use in solving problems in geometry and mechanics. For example, one problem that baffled many was the so-called ``inverse problem of tangents'': this is the problem of finding the equation of the curve  for which the tangent has specified properties. We present below two examples from mechanics.
\subsubsection{Catenary: The hanging chain}
Here the problem is to find the shape assumed by a flexible chain suspended between two points and hanging under its own weight. In 1690, Jakob Bernoulli, brother of Johann, published this problem as
challenge to the scientific world. John
Bernoulli, Gottfried Leibniz, and Christiaan Huygens each
independently solved the problem. All three solutions were
published in 1691. Burnoulli's solution used differential equations of the variable separable type. He could not complete the integration because the exponential and logarithmic functions were then unknown. (For an account of Bernoulli's approach in his own words, see \cite{bernoulli}.) 

Let the $y$-axis pass through the lowest point of the chain, let $s$ be the arc length from this point to a variable point $(x,y)$ and let $w_0$ be mass per unit length of the string (assumed constant). 

Thee portion of the chain between the lowest point and $(x,y)$ is in equilibrium under the action of three forces: the horizontal tension $T_0$ at the lowest point, the variable tension $T$ at $(x,y)$ which acts along the tangent because of the flexibility of the curve and the downward force due to the weight of the chain between these two points. Equating the horizontal and vertical components of the forces we get
$$
T\cos\theta = T_0,\quad T\sin\theta = w_0sg
$$
where $g$ is the acceleration due to gravity. Dividing we get
$$
\tan \theta  = \left(\frac{w_0g}{T_0}\right)s
$$
that is,
$$
y^\prime  = \left(\frac{w_0g}{T_0}\right)s
$$
Differentiating we get
$$
y^{\prime\prime} = \left(\frac{w_0g}{T_0}\right)\frac{ds}{dx} = \left(\frac{w_0g}{T_0}\right)\sqrt{1+(y^\prime)^2}
$$
Thus we get the differential equation of the desired curve as 
$$
y^{\prime\prime} = \left(\frac{w_0g}{T_0}\right)\sqrt{1+(y^\prime)^2}
$$
To solve this, we write $y^\prime =p$ and get
$$
p^\prime = a \sqrt{1+p^2}
$$
where $a= \frac{w_0g}{T_0}$. Separating the variables we have
$$
\frac{dp}{\sqrt{1+p^2}} = adx
$$
Solving we get 
$$ \sinh^{-1}p = ax + c
$$
Using the fact that $p=0$ when $x=0$ we get
$$
p = \sinh ax
$$
A second integrtion yields
$$
y = \frac{1}{a}\cosh ax 
$$
where we have assumed that $y=\frac{1}{a}$ when $x=0$. Thus the equation of the curved assumed by a uniform flexible chain hanging under its own weight is
$$
y = \frac{1}{a}\cosh ax 
$$
This curve is called the catenary. But the hyperbolic cosine wasn’t known by a specific
expression or name until 1761 when Lambert\footnote{Johann Heinrich Lambert (1728 – 1777) was a Swiss polymath who made important contributions to the subjects of mathematics, physics, philosophy, astronomy and map projections.} introduced the terminology and definitions.

\subsubsection{Elastica curve: The bending beam}
James Bernoulli posed the elastica problem in 1691 thus(see \cite{Levien}):

``Assuming a lamina AB of uniform thickness and width and negligible weight of its own,
supported on its lower perimeter at A, and with a weight hung from its top at B, the force
from the weight along the line BC sufficient to bend the lamina perpendicular, the curve of
the lamina follows this nature: at every point along the curve, the
product of the radius of curvature and the distance from the line BC is a constant.''
The equation of the curve of the lamina was determined as the solution of the following differential equation:
$$
\frac{dy}{dx} = \frac{x^2}{\sqrt{a^4-x^4}}
$$
But the integral cannot be evaluated in terms of elementary functions.

\subsection{Integrating factors}
The idea of using an integrating factor to solve differential equations was due to John Bernoulli, though he did not introduce the terminology. He recorded his ideas in the lecture notes he prepared in the early 1690's but published only in 1742. He used the idea to solve the equation
$$
2ydx - xdy =0.
$$
Bernoulli could not solve this equation by separation of variables because he did not know the integral of $\frac{1}{x}$. Multiplying the equation by $\frac{x}{y^2}$ and noting that
$$
d\left(\frac{x^2}{y}\right)= (2ydy-dx)\frac{x}{y^2}
$$
it was shown that the given equation is equivalent to
$$
d\left(\frac{x^2}{y}\right) =0
$$
which can be solved by direct integration to yield the solution as
$$
\frac{x^2}{y}=a,\quad \text{some constant.}
$$
Bernoulli also noted that $\frac{y^{a-1}}{x^2}$ is an integrating factor of the equation
$$
axdy - ydx =0.
$$

Another differential equation that Bernoulli solved by his method of multiplying by a factor was
$$
3y dx = x dy + y dx.
$$

Though Bernoulli introduced the concept implicitly, it was Euler\footnote{Leonhard Euler (1707 – 1783) was a Swiss mathematician, physicist, astronomer, logician and engineer who made important and influential discoveries in many branches of mathematics like infinitesimal calculus and graph theory while also making pioneering contributions to several branches such as topology and analytic number theory. He also introduced much of the modern mathematical terminology and notation, particularly for mathematical analysis, such as the notion of a mathematical function.} who developed the full theory of integrating factors in a paper published in 1741. In the paper, as an example, Euler solved the following differential equation in much the same way as it is now being solved in class room:
$$
dt + 2tz dt - t dz + t^2 dz = 0.
$$

\section{Singular solutions}
The singular solutions are not generally emphasised in the undergraduate curriculum. In this section, we define the concept, give examples and then give the historical problem that led to the formation of the concept.

In the field of differential equations, an {\em initial value problem} (also called the Cauchy problem by some authors) is the problem of finding a solution of an ordinary differential equation having specified values, that is, satisfying certain initial conditions, at a given point in the domain of the solution.

Singular solutions have been defined differently by various authors:
\begin{itemize}
\item
A singular solution $y_s(x)$ of an ordinary differential equation is a solution for which the initial value problem  fails to have a unique solution at {\em some point} on the solution. 
\item
A singular solution $y_s(x)$ of an ordinary differential equation is a solution for which the initial value problem  fails to have a unique solution at {\em every point} on the curve.
\item
The singular solution $y_s(x)$ of an ordinary differential equation is the envelope of the family of solutions.
\end{itemize}
\subsection{Example}
Consider the differential equation
\begin{equation}\label{eq1}
(y^\prime)^2 = 4y.
\end{equation}
It can be shown that a general solution of this equation can be given as 
\begin{equation}\label{eq2}
y = (x-c)^2
\end{equation}
where $c$ is an arbitrary constant. It can be seen that the following is also a solution which cannot be obtained as a particular case of the general solution by assigning any specific value to  the arbitrary constant $c$:
\begin{equation}\label{eq3}
y=0.
\end{equation}
It can be seen that the solution given by Eq.\eqref{eq2} is a singular solution in the sense of the first definition given above and that given by Eq.\eqref{eq3} in the sense of the remaining definitions.

The differential equation Eq.\eqref{eq1} still many other solutions. For example, if we choose arbitrary constants $c_1$ and $c_2$ such that $-\infty\le c_1\le c_2\le \infty$ then the following is a solution of Eq.\eqref{eq1}.
\begin{equation*} 
y  = \begin{cases} (x-c_1)^2 & \text{ for } x< c _1 \\
0 & \text{ for } c_1\le x\le c_2\\
(x-c_2)^2 & \text{ for } x>c_2
\end{cases}\\
\end{equation*}
It can be shown that this the general solution of Eq.\eqref{eq1}. But, note that it contains two arbitrary constants $c_1$ and $c_2$. 
\subsection{Discovery of singular solutions}
Singular solutions were discovered in a rather surprising
manner. Brook Taylor (of ``Taylor series" fame!), in 1715, was trying to discover the solution of the following differential equation:
$$
(1+x)^2 (y^\prime)^2 =4y^3-4y^2
$$
To solve the equation, Taylor made the substitution
$$
y = u^mv^n
$$
where $u, v$  are variables and $m,n$ are constants to be determined. Tayor chose $v=1+x^2, m= -2, n=1$ and reduced the given equation to
$$
u^2-2xu\frac{du}{dx} +v\left(\frac{du}{dx}\right)^2 =1.
$$
If we differentiate this equation we get
$$
2\frac{d^2u}{dx^2}\left(v\frac{du}{dx}-xu\right)=0.
$$
This reduces to the pair of equations:
\begin{align*}
\frac{d^2u}{dx^2} & = 0\\
v\frac{du}{dx}-xu & = 0
\end{align*}
The former gives 
$$
\frac{du}{dx}=c
$$
Substituting this in the differential equation for $u$ and solving the resulting algebraic equation for $u$ and then substituting in 
$u^mv^n$ we get the general solution of the equation as
$$
y = \frac{1+x^2}{(ax + \sqrt{1-a^2})^2}.
$$
The latter, along with the differential equation for $u$ yields the following solution
$$y=1$$
which cannot be obtained from the general solution. 
Thus we have a singular solution. This is so trivial that one could have guessed this solution!

The real significance  of the discovery singular solutions in the conceptual framework was that mathematicians had not grasped fully what is to be understood by ``solution'' of an ordinary differntial equation (see \cite{Bell} p.403).
\section{Higher order differential equations}
By the end of the seventeenth century practically all the 
known elementary methods of solving equations 
of the first order had been brought to light 
(see \cite{Ince} p.532). 

The early years of the eighteenth century are remarkable 
for a number of problems which led to differential equations 
of the second or third orders. In 1696 James Bernoulli formulated
 the isopcrimctric problem, or the problem of determining curves 
of a given perimeter which shall under given conditions,
enclose a maximum area. Five years later he published his 
solution, which
depends upon a differential equation of the third order.
Attention was now turned to trajectories in a general sense 
and in particular to
trajectories defined by the knowledge of how the curvature varies from point to
point; these gave rise to differential equations of the second order. Thus, in 1716, John Bernoulli,  discussed
an equation which would now be written in the following form:
$$
\frac{d^2 y}{dx^2} =\frac{2y}{x^2}.
$$
\section{Linear differential equations of the second order}
For simplicity, let us consider linear second order differential equations. The general form of such an equation is
\begin{equation}\label{eq8}
y^{\prime\prime}+P(x)y^\prime +Q(x)y = R(x)
\end{equation}
Equations of this kind are are important in physics, and the theory of electrical circuits. In general, the equation cannot be solved in terms of known elementary functions or even in terms of indicated integrations. However, there are methods for solving the equation in certain very special cases like when $P(x)$ and $Q(x)$ are constants or when $P(x)=ax$ and $Q(x)=bx^2$ for $a$ and $b$. Euler, as early as 1739, had come up with such methods. An important later method is to use power series expansions of $y(x)$.  

However, we have the following result regarding the existence and uniqueness of solutions of equations of the form Eq.\eqref{eq8}. It may be curious to note that Simmons closes his book (see \cite{Simmons}p.435) and signs off with a proof of this theorem!
\begin{theorem}
Let $P(x)$, $Q(x)$ and $R(x)$ be continuous functions in a closed interval $[a,b]$. If $x_0$ is any point in $[a,b]$, and and if $\alpha$ and $\beta$ be any numbers, then Eq.\eqref{eq8} has one and only one solution $y(x)$ on the on the interval $[a,b]$ such that $y(x_0)=\alpha$ and $y^\prime(x_0)=\beta$.
\end{theorem}

This theorem tells us that the solution $y(x)$ is completely determined by the values of $y(x)$ and $y^\prime(x)$ at some point $x_0$ in $[a,b]$.

The next theorem gives us a method for finding the solutions of Eq.\eqref{eq8} in the special case when $R(x)=0$.
\begin{theorem}
Let $y_1(x)$ and $y_2(x)$ be solutions of the equation
\begin{equation}\label{eq9}
y^{\prime\prime}+P(x)y^\prime +Q(x)y =0
\end{equation}
on the interval $[a,b]$. Then 
$$
c_1y_1(x)+c_2y_2(x)
$$
is a general solution of the equation in the sense that every solution of the equation in $[a,b]$ can be obtained by a suitable choice of $c_1$ and $c_2$.
\end{theorem}
By Theorem 1, a solution of Eq.\eqref{eq9} is completely determined by the values of $y(x_0)$ and $y^\prime(x_0)$ at some point $x_0$ in $[a,b]$. So, we have to show that, given arbitrary $\alpha$ and $\beta$, we have a solution of the form $y(x)=c_1y_1(x)+c_2y_2(x)$ with $y(x_0)=\alpha$ and $y^\prime(x_0)=\beta$. Since the solution is unique, and since $c_1y_1(x)+c_2y_2(x)$ is a solution, it is enough to prove that we can find $c_1$ and $c_2$ satsfying the following conditions:
\begin{align*}
c_1y_1(x_0)+c_2y_2(x_0) & = \alpha\\
c_1y_1^\prime(x_0)+c_2y_2^\prime(x_0)& = \beta
\end{align*}
This is possible if and only if 
$$
\begin{vmatrix} y_1(x_0) & y_2(x_0) \\ y_1^\prime (x_0) & y_2^\prime(x_0) \end{vmatrix} \ne 0.
$$
This must be true for all $x_0$ in $[a,b]$.
This assured by the next theorem where we have used the function
$$
W(y_1(x), y_2(x))=\begin{vmatrix} y_1(x) & y_2(x) \\ y_1^\prime (x) & y_2^\prime(x) \end{vmatrix} 
$$
called the Wronskian of $y_1(x)$ and $y_2(x)$.
\begin{theorem}
Let $y_1(x)$ and $y(x_2)$ be solutions of Eq.\eqref{eq8} on $[a,b]$, then the Wronskian $W(y_1(x), y_2(x))$ is either identically zero, nor never zero in $[a,b]$.
\end{theorem}

A proof of this is instructive. We begin by noting  that 
$$
\frac{dW}{dx} = y_1(x)y_2^{\prime\prime}(x) - y_2(x)y_1^{\prime\prime}(x).
$$
Now observing that $y_1(x)$ and $y_2(x)$ are solutions of Eq.\eqref{eq8}, we have
$$
\frac{dW}{dx}=-PW.
$$
Solving this differential equation, we have
$$
W=ce^{-\int P\, dx}.
$$
Since the exponential is never zero, the result follows.
\section{Linear differential equations of higher orders}
\subsection{Homogeneous equations with constant coefficients}
The general treatment of homogeneous linear differential equations with constant coefficients could be said to have been inaugurated by Euler with a letter written to John Burnoulli on 15 September 1739. Euler discovered all the methods for solving such equations that are currently taught in undergraduate classes.
\begin{equation}\label{eq7}
0=a_0y +a_1\frac{dy}{dx}+a_2\frac{d^2y}{dx^2} + \cdot + a_n\frac{d^n y}{dx^n}
\end{equation}

Let us follow Euler's method as explained by Ince (see \cite{Ince} p.585). 
\begin{enumerate}
\item
If $y=u$ is a solution of Eq.\eqref{eq7} the $y= cu$ is a also a solution of the same equation where $c $ is any constant.
\item
If we can obtain $n$ particular solutions $y=y_1$, $y=y_2$, $\ldots$, $y=y_n$, the complete or general solution will be 
$$
y=c_1y_1+\cdot+y_ny_n
$$
where $c_1,\ldots,c_n$ are constants.
\item
If $z=\alpha$ is a solution of the equation
\begin{equation}\label{eq8}
a_0+a_1z +\cdots +a_nz^n=0
\end{equation}
then $y=e^{\alpha x}$ will satisfy Eq.\eqref{eq7}.
\item
There are as many particular solutions of this form as there are real factors of the form $z-\alpha$ in 
$$
a_0+a_1z +\cdots +a_nz^n.
$$
\item
If there is a multiple factor of the form $(z-\alpha)^k$, then
using the substitution $y=e^{\alpha x}u$, a solution involving $k$ constants can be found:
$$
y=e^{\alpha x}(c_1+c_2x+\cdots + c_kx^k).
$$
\item
When a pair of complex factors arise, they are united in a real quadratic factor of the form $p-qz+rz^2$ which corresponds to the differential equation
$$
0=py-2z\sqrt{pr}\cos \phi \frac{dy}{dx}+r\frac{d^2y}{dx^2}, \text{ where } \cos\phi = \frac{q}{2\sqrt{pr}}.
$$
\item
The transformation 
$y=e^{\sqrt{pr} x \cos\phi }u$
reduces the equation to an equation of the form
$$
\frac{d^2y}{dx^2}+Ay=0.
$$
A method for solving this equation had already been developed by Euler.
\item
The case of repeated quadratic factors was then dealt with and the discussion of the homogeneous linear equation with constant coefficients was complete.
\end{enumerate}

\subsection{Non-homogeneous equations with constant coefficients}
Euler also considered non-homogeneous linear differential equations with constant coefficients:
$$
X=a_0y +a_1\frac{dy}{dx}+a_2\frac{d^2y}{dx^2} + \cdot + a_n\frac{d^n y}{dx^n}
$$
Here $X$ is a function of $x$. The method  adopted was that of a successive
reduction of the order of the equation by the aid of integrating factors of the $e^{\alpha x}$.

We illustrate Euler's method by considering a second order equation, say,
$$
\frac{d^2y}{dx^2}+ky = X.
$$
We choose $\alpha$ such that the integral $\int e^{\alpha}x \, X\, dx$ has the following form:
$$
\int e^{\alpha}x \, X = e^\alpha\left(Ay+B\frac{dy}{dx}\right),
$$
where $A$ and $B$ are constants. To find $\alpha$, we differentiate this relation to get
\begin{align*}
e^{\alpha x} \, X & = e^{\alpha x}\left(A\frac{dy}{dx} + B\frac{d^2 y}{dx^2}\right)
+\alpha e^{\alpha x}\left(Ay + B\frac{dy}{dx}\right)\\
& = e^{\alpha x} \left( \alpha Ay + (A + \alpha B)\frac{dy}{dx} + B \frac{d^2 y}{dx^2}\right)\\
& = e^{\alpha x}\left(\frac{d^2 y}{dx^2} + ky\right)
\end{align*}
Equating the coefficients, we get
\begin{align*}
\alpha A & = k\\
A+\alpha B & = 0\\
B & =1
\end{align*}
Solving these equations we get
$$
\alpha = \sqrt{-k},\quad A = -k, \quad B = 1.
$$
Thus we have
$$
e^{-\alpha x}\int e^{\alpha x}X \, dx = -k y + \frac{dy}{dx}.
$$
$y$ can be expressed as a function of $x$ by repeating the procedure.

\subsection{Linear differential equations with variable coefficients}
Historically, it is interesting to observe that Euler studied differential equations of the following form before much before he considered equations with constant coefficients:
$$
0=a_0y +a_1x\frac{dy}{dx}+a_2x^2\frac{d^2y}{dx^2} + \cdot + a_nx^n\frac{d^n y}{dx^n}
$$
Euler did not develop a general method for solving such equations. Instead his approach was to find an ingenious substitution which reduced the equation to a similar equation of order  $n-1$.

Euler first multiplied the equation by $x^p$ and made use of the substitution
$$
z=\frac{1}{p+1}\frac{d}{dx}\left(x^{p+1}\right).
$$
The constant $p$ is appropriately chosen such that the resulting equation is a equation of order $n-1$. This process of reduction was then repeated as often as necessary.
\section{The operator $D$}
In describing methods for finding solutions of linear differential equations, we use the operator $D$ and manipulate it blindly without pausing to think what it really means. We sometimes say that $D$ represents $\frac{d}{dx}$ and state that 
$$
Dy=\frac{dy}{dx}.
$$
We also say that 
\begin{align*}
D^2 y & = \frac{d^2y}{dx^2}\\
D^3 y & = \frac{d^3y}{dx^3}\\
\vdots & 
\end{align*}
 As a typical example, let us consider how we solve the following differential equation is solved:
$$
\frac{d^2y}{dx^2}-3\frac{dy}{dx} + 2y = x +1
$$
We rewrite thw equation using $D$ as follows:
$$
(D^2-3D+2)y = x + 1
$$
The auxiliary equation
$$
D^2-3D+2=0
$$
is then set up and its solutions are found as $D=1,2$. From this the complementary function is stated as
$$
\text{CF } = c_1 e^{x}+c_2e^{2x}
$$
where $c_1$ and $c_2$ are arbitrary constants. A particular integral is determined in a more bizarre way.
\begin{align*}
\text{PI } 
& = \frac{1}{D^2-3D+2}(x+1)\\
& = \frac{1}{2}\frac{1}{1+\frac{D^2-3D}{2}} (x+1)\\
& = \frac{1}{2}\left( 1 - \left(\frac{D^2-3D}{2}\right) + \left(\frac{D^2-3D}{2}\right)^2 - \left(\frac{D^2-3D}{2}\right)^3 + \cdots \right)(x+1)\\
& = \frac{1}{2}\left( (x+1) - \frac{1}{2}(D^2(x+1) -3D(x+1)) + 0 + \cdots\right)\\
& = \frac{1}{2}\left( (x+1) - \frac{1}{2}(0 - 3)\right)\\
& = \frac{1}{2}\left(x+\frac{5}{2}\right)
\end{align*}
It is now claimed that the complete solution of the equation is 
\begin{align*}
y & = \text{ CF } + \text{ PI }\\
& = c_1 e^{x}+c_2e^{2x} + \frac{1}{2}(x+\frac{5}{2})
\end{align*}
OK. Fine. We do have a solution of the given differential equation. But, what is happening? 

Let $A$ be the set of all real valued differentiable functions defined over some domain and $B$ be the set of all derivatives of such differentiable functions. Then $D$ is a mapping from $A$ to $B$ defined by 
$$
D : f(x) \mapsto \frac{d}{dx}f(x).
$$
In the elementary theory of differential equations, we consider a set smaller than $A$. We consider the set $S$ of all analytic functions where by analytic we mean that the function is
infinitely differentiable which in turn means that the function possesses derivatives of every order and then consider $D$ as a mapping from $S$ to itself. 

It can be easily verified that $S$ has the structure of a real vector space and that $D$ is a linear operator on $S$. Since $D$ 
is a linear operator on $S$, we can consider the composition $D\circ D$ of $D$ with itself. This composition is denoted by $D^2$. Thus $D^2$ is the not the square of $D$ in the sense of ordinary multiplication. The mappings $D^3$, $D^4$, $\ldots$ are defined in a similar way.
\begin{align*}
D^2 (f(x)) 
& = (D\circ D)(f(x)) \\
& = D(D(f(x)))\\
& = \frac{d}{dx}\left(\frac{d}{dx}f(x)\right)\\
& = \frac{d^2}{dx^2}f(x)
\end{align*}
It is well known that the set of linear operators on a real vector space itself has the structure of a real vector space and so we can consider linear combinations of the operators $D, D^2, \ldots$ and form operators of the form
\begin{equation}\label{Eq4}
a_0 + a_1D+a_2D^2 + \cdots a_nD^n.
\end{equation}

The operator $D$ is a many-one mapping and it has no inverse. However, given any $f(x)$ in $S$, we can define $D^{-1}(f(x))$ by
$$
D^{-1}(f(x)) = \{g(x)| D(g(x)) = f(x)\}.
$$
Using the integral notation, this can be expressed in the form
$$
D^{-1}(f(x)) =\int f(x)\, dx
$$
This idea can be extended to operators of the more general form given in \eqref{Eq4} also. This justifies the use of the notation
$(D^2-3D+2)^{-1}$, which may be expressed as $\frac{1}{D^2-3D+2}$.

Some advanced analysis is required to see why the operator $\frac{1}{1+\frac{D^2-3D}{2}}$ may expanded in powers of  $\frac{D^2-3D}{2}$. 
That some restrictions are required for the validity of this expansion can be seen by considering the following problem.

Let us find a particular integral of the following equation by the method indicated above:
$$
\frac{dy}{dx}+y= \frac{1}{x}
$$
We have
\begin{align*}
\text{PI }
& = \frac{1}{1+D}\frac{1}{x}\\
& = (1-D+D^2-D^3+\cdots)\frac{1}{x}\\
& = \frac{1}{x} + \frac{1}{x^2}+\frac{2!}{x^3}+\frac{3!}{x^4}+\cdots\\
& = \sum_{n=0}^\infty \frac{n!}{x^{n+1}}
\end{align*}
But the last series is not convergent for any value of $x$ and so it does not even represent a function, let alone a solution of the differential equation.

\section{End of an era}
The period of  discovery of general methods for solving ordinary  differential equations ended by 1775, a hundred years after Leibniz inaugurated the integral sign. For many problems the formal methods were not sufficient. Solutions with special properties were required, and thus, criteria guaranteeing the existence of such solutions became increasingly important. Boundary value problems led to ordinary differential equations, such as Bessel's equation, that prompted the study of Laguerre, Legendre, and Hermite polynomials. The study of these and other functions that are solutions of equations of hypergeometric type led in turn to modern numerical methods.

\section{Solutions in series}
It was Newton who first obtained a power series solution of a differential equation. He indicated it as general method and illustrated the idea by solving a few problems. Though Newton wrote about these ideas in 1671, it was published in print only in 1736. 

The class of elementary functions consists of algebraic functions, the trigonometric,  inverse trigonometric,  exponential and logarithmic functions and all others that can be constructed from these by adding, subtracting, multiplying, dividing, or forming function of a function. Beyond the elementary functions lie the {\em higher transcendental functions}, also called {\em special functions}. Since the beginning of the eighteenth century, many hundreds of special functions have been considered sufficiently interesting or important to merit some degree of study. Most of them are now completely forgotten. A few have survived because of their applications and intrinsic values. 

A large class of special functions arises as solutions of second order linear differential equations. The method of power series is used to obtain solutions of such equations and and the resulting solutions are taken as the definitions of special functions. An understanding of the concepts of convergence of series is essential for a proper appreciation of the definitions. 

We illustrate the general procedure by a typical example. We take a very formal approach ignoring questions of convergence. However, it is advisable to have a look at the theorems which guarantee the validity of these process (for example, see Theorem 29A in \cite{Simmons}).

\subsection{Gauss's hypergeometric equation}
The following differential equation is known as Gauss's\footnote{Johann Carl Friedrich Gauss (1777 –  1855) was a German mathematician who contributed significantly to many fields, including number theory, algebra, statistics, analysis, differential geometry, geodesy, geophysics, mechanics, electrostatics, magnetic fields, astronomy, matrix theory, and optics. Sometimes referred to as the {\em Princeps mathematicorum} (Latin, ``the foremost of mathematicians") and ``greatest mathematician since antiquity", Gauss had an exceptional influence in many fields of mathematics and science and is ranked as one of history's most influential mathematicians.}
hypergeometric equation:
\begin{equation}\label{eq6}
x(1-x)y^{\prime\prime} +[c-(a+b+1)x]y^\prime - ab y=0
\end{equation}
The coefficients are chosen in such a way that the solution has a nice form.

By theorem (Theorem 29A in \cite{Simmons}), the equation has a solution of the following form:
$$
y=a_0+a_1x+\cdots+a_nx^n+\cdots = \sum_{n=0}^\infty a_nx^n
$$
where $a_0$ is a nonzero constant. We have
\begin{align*}
y^\prime & = \sum_{n=1}^\infty na_nx^{n-1}\\
y^{\prime\prime} & = \sum_{n=2}^\infty n(n-1)a_{n-2}x^{n-2}
\end{align*}
Substituting these in Eq,\eqref{eq6}, and equating the coefficients of $x^n$ to $0$ we get
$$
a_{n+1}=\frac{(a+n)(b+n)}{(n+1)(c+n)}a_n.
$$
Setting $a_0=1$, we get one solution as
$$
y=1+\frac{ab}{1\cdot c}x+\frac{a(a+1)b(b+1)}{n(n+1)c(c+1)}x^2+\cdots
$$
This series is the {\em hypergeometric series} is denoted by $F(a,b,c,x)$. It can be shown that when $c$ is not $0$ or  a negative integer, $F(a,b,c,x)$ is analytic function in the interval $-1<x<1$. It is then called the hypergeometric function. If $c$ is not a positive integer, thee hypergeometric function has a second linearly independent solution given by
$$
y= x^{1-c}F(a-c+1, b-c+1, 2-c, x).
$$ 
Thus the general solution of the hypergeometric equation is given by
$$
y = c_1F(a,b,c,x) + c_2 x^{1-c}F(a-c+1, b-c+1, 2-c, x).
$$
 
It is instructive to note that the elementary transcendental functions can be expressed in terms of the hypergeometic function. Thus the hypergeometric function ``unifies'' all elementary transcendental functions.
\renewcommand{\arraystretch}{1.5}
\begin{table}
\begin{center}
\begin{tabular}{ll}
\hline
Function & In terms of hypergeometric function\\
\hline
$(1+x)^p $& $F(-p,b,b,-x)$\\
$\log (1+x)$ & $ xF(1,1,2,-x)$\\
$\sin^{-1}x$ & $xF\left(\frac{1}{2},\frac{1}{2}, \frac{3}{2}, x^2\right)$ \\
$\tan^{-1} x$ & $xF\left(\frac{1}{2},\frac{1}{2}, 1, -x^2\right)$ \\
$e^x$ & $\lim_{b\rightarrow \infty}F\left(a,b,a, \frac{x}{b}\right)$\\
$\sin x$ & $x\left[ \lim_{a\rightarrow \infty} F\left(a,a,\frac{3}{2}, \frac{-x^2}{4a^2}\right)\right]$\\
$\cos x$ & $\lim_{a\rightarrow \infty }F\left(a,a,\frac{1}{2}, \frac{-x^2}{4a^2}\right)$\\
\hline
\end{tabular}
\caption{Elementary transcenental functions in terms of the hypergeometric functions}
\end{center}
\end{table}
\section{Nonlinear differential equations}
So far we have confined ourselves with linear differential equations mostly because of the fact that such equations are more
amenable to finding solutions, even though this is rarely possible.  
In this section we have a cursory look at the theory of nonlinear differential equations. In this theory, no attempts are made to obtain solutions in the traditional sense instead the efforts are to obtain qualitative information about the general behavior of solutions. The qualitative theory of nonlinear equations was founded by Poincare around 1880 in connection with his work on celestial mechanics. 

Our attempt here is only to give a flavor of the type of equations considered, the type of questions asked and the nature of answers obtained for such questions. 

An important class of nonlinear equations consists of systems of equations of the following form:
\begin{align}
\frac{dx}{dt} & = F(x,y) \label{eq10}\\
\frac{dy}{dt}& = G(x,y)\label{eq11}
\end{align}
A sysytm of this kind, in which the independent variable $t$ does not appear in the functions $F(x,y)$ and $G(x,y)$ are called {\em autonomous} systems. One of the well-known system of this type is Volterra's prey-predator equations (published  by Vito Volterra\footnote{Vito Volterra (1860 – 1940) was an Italian mathematician and physicist, known for his contributions to mathematical biology and integral equations, being one of the founders of functional analysis.} in 1926):
\begin{align*}
\frac{dx}{dt} & = ax- bxy\\
\frac{dy}{dt}& = -cy + dxy
\end{align*}
These equations are used to describe the dynamics of biological systems in which two species interact, one as a predator and the other as prey.\footnote{The equations are derived based on the following assumptions:
\begin{itemize}
\item
    The prey population finds ample food at all times.
\item
    The food supply of the predator population depends entirely on the size of the prey population.
\item
    The rate of change of population is proportional to its size.
\item
    During the process, the environment does not change in favour of one species and genetic adaptation is inconsequential.
\item
    Predators have limitless appetite.
\end{itemize}}

Let $x=x(t)$, $y=y(t)$ be a solution of the system Eq.\eqref{eq10}-\eqref{eq11}. Then as $t$ varies the point $(x(t), y(t))$ traces out a curve in the $xy$-plane. Such a curve is called a {\em path} of the system. At most one path passes through each point in the plane. The points $(x_0,y_0)$ where $F(x_0,y_0)=0$ and $G(x_0,y_0)=0$ are special and they are called the {\em critical points} of the system.

In regard to an autonomous sytem like Eq.\eqref{eq10}-\eqref{eq11} the following questions are posed:
\begin{enumerate}
\item
What are the critical points?
\item
How are the paths near the critical points arranged?
\item
Does a point near a critical point remains near or wanders off into another part of the plane as $t$ increases (the stabilty or instability of critical points)? 
\item
Are there closed paths enclosing the critical points? (Such path correspond to periodic solutions.)
\end{enumerate}

By assuming that $(0,0)$ is a critical point and expanding $F(x,y)$ and $G(x,y)$ as power series in $x$ and $y$ and then retaining only 
terms of the first degree the general autonomous system represented by Eq.\eqref{eq10}-\eqref{eq11} can be approximated by a system of the following form:
\begin{align*}
\frac{dx}{dt} & = a_1x+b_1y\\
\frac{dy}{dt} & = a_2 x + b_2 y
\end{align*}
It will be assumed that
$$
a_1b_2-a_2b_1\ne 0.
$$

It can be shown that this system has a solution of the form
\begin{align*}
x& = Ae^{mt}\\
y & = Be^{mt}
\end{align*}
where $m$ is a root of the quadratic equation
$$
m^2 - (a_1+b_2)m +(a_1b_2-a_2b_1) =0.
$$
Let $m_1$ and $m_2$ be the roots of this quadratic equations. The nature of the critical point $(0,0)$ is determined by the nature of the numbers $m_1$ and $m_2$. The following cases arise:
\begin{enumerate}
\item
$m_1$ and $m_2$ are real, distinct and of the same sign ({\em node}).
\item
$m_1$ and $m_2$ are real, distinct and of opposite signs ({\em saddle point}).
\item
$m_1$ and $m_2$ are conjugate complex, but not pure imaginary ({\em spiral}).
\item
$m_1$ and $m_2$ are real and equal ({\em node}).
\item
$m_1$ and $m_2$ are pure imaginary ({\em center}). 
\end{enumerate}
The {\em phase portrait} of an autonomous system is a diagram giving an overall picture of the paths. By constructing and analysing the phase portrait one can study the stability properties of the paths.
\addcontentsline{toc}{section}{References}


\begin{thebibliography}{99}
\bibitem{Ince}
E. L. Ince, {\em Ordinary Differential Equations}, Publisher: Not available, 1920. Full text is available  at {\tt https://archive.org/details/ ordinarydifferen029666mbp}.
\bibitem{Sasser}
John E. Sasser, ``History of ordinary differential equations:  
The first hundred years''. Available  at {\tt http://www2.fiu.edu/yuasun/ ODE\_History.pdf}.
\bibitem{Levien}
Raph Levien, ``The elastica: a mathematical history'', Technical Report No. UCB/EECS-2008-103, August 23, 2008. Available  at {\tt 
http://www.eecs.berkeley.edu/Pubs/TechRpts/2008/ EECS-2008-103.html}
\bibitem{Bell}
E. T. Bell, {\em The development of Mathematics}, McGraw Hill Book Company, 1945.
\bibitem{Phaser}
Huseyin Kocak, ``Newton's first differential equation''. Available at
{\tt http://www.phaser.com/modules/historic/newton/ isaac\_newton\_ode.pdf} 
\bibitem{Newton}
Issac Newton,  [1736].
{\em The Method of Fluxions and Infinite Series; with its Applications to the Geometry of Curve-lines}
by the Inventor Sir Isaac Newton,
Kt., Late President of the Royal Society. Translated from the Author’s Latin Original not yet made publick. To which is subjoin’d, A perpetual Comment upon the whole Work, Consisting of
Annotations, Illustrations, and Supplements, In Order to make this Treatise A Compleat Institution for the use of Learners. By John Colson. London: Printed by Henry Woodfall; And Sold by John Nourse, at the Lamb without Temple-Bar. M.DCCXXVI. Available at 
{\tt https://archive.org/details/methodoffluxions00newt}
\bibitem{Simmons}
George F. Simmons, {\em Differential Equations with Applications and Historical Notes}, 2nd Edition, McGraw Hill Education, July 2017.
\bibitem{bernoulli}
Johann Bernoulli 
(Translated by William A. Ferguson, Jr.), ``Lectures on 
The Integral
Calculus'', 21st Century Science \& Technology, 2004. Available at {\tt http://21sci-tech.com/translations/Bernoulli.pdf}.
\bibitem{Coddington}
Earl A. Coddington and Norman Levinson, {\em Theory of Ordinar Differential equation}, McGraw Hill, 1955.

\end{thebibliography}
\end{document}